\documentclass[11pt,twoside]{article}
\usepackage{a4wide}

\usepackage{amssymb,amsmath}

\providecommand{\proof}{{\bf Proof: }}

\newtheorem{lemma}{LEMMA}
\newtheorem{theorem}{THEOREM}
\newtheorem{remark}{REMARK}
\newtheorem{example}{EXAMPLE}

\providecommand{\eps}{\varepsilon}

\providecommand{\N}{\mathbb{N}}
\providecommand{\R}{\mathbb{R}}
\providecommand{\AAb}{{\mathbf A }}
\providecommand{\bb}{{\mathbf b }}
\providecommand{\fT}{f_\mathcal{T}}
\providecommand{\uT}{u_\mathcal{T}}
\providecommand{\vT}{v_\mathcal{T}}
\providecommand{\wT}{w_\mathcal{T}}
\providecommand{\xT}{x_\mathcal{T}}
\providecommand{\zT}{z_\mathcal{T}}
\providecommand{\uTi}{u_{\mathcal{T}i}}
\providecommand{\vTi}{v_{\mathcal{T}i}}
\providecommand{\wTi}{w_{\mathcal{T}i}}
\providecommand{\zTi}{z_{\mathcal{T}i}}
\providecommand{\uTj}{u_{\mathcal{T}j}}
\providecommand{\vTj}{v_{\mathcal{T}j}}
\providecommand{\wTj}{w_{\mathcal{T}j}}
\providecommand{\zTj}{z_{\mathcal{T}j}}
\providecommand{\zetaTij}{\zeta_{\mathcal{T}ij}}
\providecommand{\II}{{\mathbf I }}
\providecommand{\IT}{I_\mathcal{T}}
\providecommand{\LT}{L_\mathcal{T}}
\providecommand{\VT}{V_\mathcal{T}}
\providecommand{\close}{\nolinebreak\mbox{$\blacktriangleleft$}}
\providecommand{\inte}{\mathop{\rm int\,}}
\providecommand{\measd}[1]{\mbox{meas}_d\left(#1\right)}
\providecommand{\sumi}{\sum_{i\in\Lambda}}
\providecommand{\sumj}{\sum_{j\in\Lambda_i}}
\newcounter{property}
\newcommand{\refsteppropl}[1]{\refstepcounter{property}\label{#1}}
\renewcommand{\theproperty}{P\arabic{property}}
\title{%
\emph{A posteriori} estimates for errors of functionals on finite volume approximations
to solutions of elliptic boundary value problems}
\author{Lutz Angermann\thanks{Technische Universit\"at Clausthal,
Institut f\"ur Mathematik,
Erzstr. 1, D--38678 Clausthal-Zellerfeld,
Federal Republic of Germany,
e-mail: lutz.angermann@tu-clausthal.de%
}}
\date{}
\begin{document}
\bibliographystyle{alpha}
\maketitle
{\small
\emph{Abstract:}
This article describes the extension of recent methods
for \emph{a posteriori} error estimation such as dual-weighted
residual methods to node-centered finite volume discretizations
of second order elliptic boundary value problems
including upwind discretizations.
It is shown how different sources of errors,
in particular modeling errors and discretization errors,
can be estimated with respect to a user-defined output functional.}

\vspace{.8ex}
{\small
\emph{Keywords:}
Finite volume methods,
a posteriori error estimates,
DWR method,
goal-oriented estimation}

\vspace{.8ex}
{\small
\emph{2010 Mathematics Subject Classification:}
65\,N\,08, % Finite volume methods
65\,N\,15, % Partial differential equations, boundary value problems:
%               % Error bounds
65\,N\,30, % Finite elements, Rayleigh-Ritz and Galerkin methods,
%               % finite methods
65\,N\,50% % Mesh generation and refinement
}

\section{Introduction}

In many areas of practical interest, e.g.\ computational fluid dynamics or image
reconstruction, the computations are complicated and expensive, effectively
limiting the achievable precision. In order to overcome these problems, adaptive
finite element approaches are in use since several decades
(see, e.g., \cite{Babuska:78}). For instance, in the so-called $h$-adaptive
methods the computational meshes are refined locally so that the mesh captures the
variation of the solution while remaining coarse elsewhere. It has been shown
that such approaches are computationally much more efficient than uniform
meshes.
In recent years, there has been considerable progress in applying these techniques
to more involved questions such as the \emph{a posteriori} error estimation
of values of nonlinear functionals of interest (\emph{goal-oriented} estimation,
see, e.g., \cite{Becker:01b}, \cite{Bangerth:03}, \cite{Rannacher:05})
or the (additional) \emph{a posteriori} estimation
of modeling errors (see, e.g., \cite{Oden:00}, \cite{Oden:01}, \cite{Braack:03}).

The present paper describes the extension of recent techniques
for obtaining \emph{a posteriori} error estimates
for modeling and discretization errors
to nonlinear second-order elliptic PDEs
which are discretized by means of node-centered finite volume schemes
including stabilization mechanisms of upwind type.
Finite volume methods are attractive methods in selected areas
of application, and therefore it is a natural requirement to develop
analogous methods of error control for FVM. However, since finite volume methods
suffer, in general, from the so-called property of Galerkin-orthogonality,
special attention is to be paid to the treatment of the resulting defect term.
It is shown that the extension of
the dual-weighted a posteriori error estimates
to finite volume discretizations is possible in a reasonable way.
Furthermore, the latter approach is interesting because of the fact
that different sources of errors
(i.e.\ not only discretization errors but, for example, also modeling errors)
can be estimated with respect to a rather arbitrary user-defined output functional.
For instance, in the field of inverse problems, the Tikhonov functionals
can serve as typical output functionals (see, e.g., \cite{Beilina:10c}).

Here we will mainly deal with Voronoi and Donald finite volume partitions
on simplicial primary partitions of the computational domain;
however the ideas can be extended to more general primary partitions,
in particular quadrilateral or hexahedral partitions
(cf., e.g., \cite[Sect.\ 4.2]{Angermann:06a}).

We consider the following boundary value problem with respect to
the unknown function $u:\;\Omega\to\R$:
\begin{equation}\label{genbvp}
\left\{\begin{array}{rcl}
-\nabla\cdot (\AAb(\cdot,u)\nabla u) + \bb(\cdot,u)\cdot \nabla u + c(\cdot,u)u
&=& f \quad \mbox{in } \Omega,\\
u &=& 0 \quad \mbox{on } \Gamma,
\end{array}\right.
\end{equation}
where $\Omega\subset\R^d,$ $d\in\{2,3\},$ is a bounded polygonal or polyhedral domain
with a Lipschitzian boundary $\Gamma,$
and the data in (\ref{genbvp}) are sufficiently smooth:
$$
\AAb:\;\Omega\times\R\to\R^{d,d},
\quad
\bb:\;\Omega\times\R\to\R^d,
\quad
c:\;\Omega\times\R\to\R,
\quad
f:\;\Omega\to\R.
$$
Equations of such type may occur in various areas of science, for example
in the mathematical description of filtration processes
in nonhomogeneous media.

Using the formal notation
\begin{eqnarray}
(w,v)&:=&\int_\Omega wv\,dx,\nonumber\\
(\nabla w,\nabla v)&:=&\int_\Omega \nabla w\cdot\nabla v\,dx,\nonumber\\
b(w;v)&:=&\frac{1}{2}\left[(\bb(\cdot,w)\cdot\nabla w,v)
-(w,\bb(\cdot,w)\cdot\nabla v)\right],\label{bcont}\\
d(w;v)&:=& (c(\cdot,w)w,v)-\frac{1}{2}\left(w\nabla\cdot \bb(\cdot,w),v\right),\nonumber\\
a(w;v)&:=& (\AAb(\cdot,w)\nabla w,\nabla v)+b(w;v)+d(w;v)\label{asplit},
\end{eqnarray}
and $\langle f,v\rangle:=(f,v),$
the variational formulation of the problem (\ref{genbvp})
in the space $V:=H^1_0(\Omega)$ reads as follows:
\begin{quote}
Find $u\in V$ such that
\end{quote}
\begin{equation}\label{genvar}
\forall v\in V:\quad a(u;v)=\langle f,v\rangle.
\end{equation}
Regarding results for the existence, uniqueness and regularity
of solutions of (\ref{genbvp}) or (\ref{genvar}),
there is a wide literature both of relatively general nature
(see, e.g., \cite[Ch.\ 2]{Boehmer:08} for a short survey) as well as for
more specialized equations (see, e.g., \cite{Antontsev:06}).

\section{The finite volume scheme}\label{s:fv_deriv}

Finite volume methods are attractive discretization methods
for partial differential equations of first or second order in conservative form
since they adequately transfer the conservation law, which is
expressed by the differential equation, to the discrete level.
At the same time, due to their proximity to finite difference methods,
they are relatively easy to implement even in the nonlinear situation.
However, a drawback of many finite volume methods is
that there is no $p$-hierarchy as in finite element methods,
therefore the order of accuracy (related to the grid size) is relatively low.
Nevertheless finite volume methods find wide applications in the
computational practice.
A certain degree of compensation can be achieved by the application
of adaptive techniques based on a posteriori error estimates,
as discussed in the subsequent section.

In this section we concentrate on node-centered finite volume methods
for the discretization of problem (\ref{genbvp}).

\subsection{The case of Voronoi diagrams and scalar diffusion coefficients}

Let us consider a family of Voronoi diagrams such that
their straight-line duals are Delaunay triangulations of $\Omega$
consisting of self-centered simplices.
Here a simplex $T$ is called \emph{self-centered} if its circumcentre
lies in the interior of $T$ or on the boundary $\partial T.$

Denote by
$\overline\Lambda\subset\N$ the index set of all vertices $x_i$
of a particular triangulation $\mathcal{T}$
and by $\Lambda\subset\overline\Lambda$
the index set of all vertices lying in $\Omega.$

\begin{figure}
\begin{center}
\setlength{\unitlength}{1mm}
\begin{picture}(120,65)(-60,-32)
% edges of triangles
\thicklines
\put(30,0){\line(-1,1){30}}
\put(0,30){\line(-3,-2){30}}
\put(-30,10){\line(1,-4){10}}
\put(-20,-30){\line(4,1){40}}
\put(20,-20){\line(1,2){10}}
\put(0,0){\line(1,0){30}}
\put(0,0){\line(-3,1){30}}
\put(-20,-30){\line(2,3){20}}
\put(20,-20){\line(-1,1){20}}
\put(0,0){\line(0,1){30}}
\thinlines
% boundary of Dirichlet region
\put(15,-5){\line(0,1){20}}
\put(15,15){\line(-1,0){26.6666}}
\put(-11.6666,15){\line(-1,-3){7.8787}}
\put(-19.5454,-8.6363){\line(3,-2){18.5454}}
\put(-1,-21){\line(1,1){16}}
% nodes
\put(0,0){\circle*{0.8}}
\put(30,0){\circle*{0.8}}
\put(0,30){\circle*{0.8}}
\put(-30,10){\circle*{0.8}}
\put(-20,-30){\circle*{0.8}}
\put(20,-20){\circle*{0.8}}
\put(-5.5,3){$x_i$}
\put(30,1.7){$x_j$}
\put(9.4,5){$\Gamma_{ij}$}
\put(30,-30){\line(1,0){4}}
\put(35,-31){\rm boundary of $\Omega_i$}
\end{picture}
\end{center}
\caption{Configuration for the Voronoi-type discretization ($d=2$)}
\label{fig:config}
\end{figure}
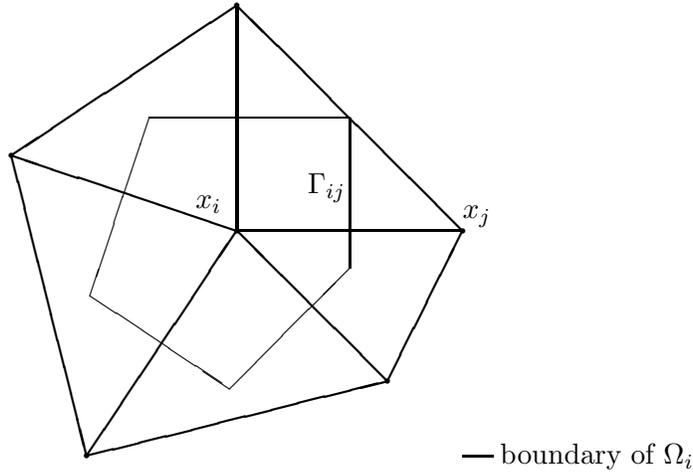

In more detail, let
\begin{eqnarray*}
\Omega_i:=\Omega_i^V&:=&\{x\in\Omega:\,\|x-x_i\|<\|x-x_i\|
\ \forall j\in\overline\Lambda\setminus\{i\}\},\ i\in\overline\Lambda,\\
&&\mbox{where }\|\cdot\| \mbox{ denotes the Euclidean norm in }\R^d,\\
m_i&:=&\measd{\Omega_i},\\
&&\mbox{where }\measd{\cdot} \mbox{ denotes the $d$-dimensional volume,}\\
\Gamma_{ij}&:=&\partial\Omega_i\cap\partial\Omega_j,
\ \Gamma_{ij}^T:=\Gamma_{ij}\cap T,
\ i\in\Lambda,j\in\overline\Lambda\setminus\{i\},\ T\in\mathcal{T},\\
m_{ij}&:=&\mbox{meas}_{d-1}\left(\Gamma_{ij}\right),
\ m_{ij}^T:=\mbox{meas}_{d-1}\left(\Gamma_{ij}^T\right),\\
d_{ij}&:=&\|x_i-x_j\|,\\
\Lambda_i&:=&\{j\in\overline\Lambda\setminus\{i\}:
\ m_{ij}\ne 0\},\\
\Lambda_T&:=&\{i\in\overline\Lambda:
\ x_i\in\partial T\},\\
h&:=&\max_{T\in\mathcal{T}}h_T,\quad \mbox{where } h_T:=\mbox{diam}\,T.
\end{eqnarray*}
The finite volume solution will be interpolated in the discrete space
$$
\VT:=\left\{v\in V:\;\left(\forall T\in\mathcal{T}:\;
v\left|_T\right.\in\mathcal{P}_1(T)\right)\right\},
$$
where $\mathcal{P}_1(T)$ is the set of all first degree polynomials on $T.$
We introduce a so called \emph{lumping operator}
$$
\LT:\ C(\overline\Omega)\to L_\infty(\Omega)
\quad\mbox{acting as}\quad
\LT v:=\sum_{i\in\overline\Lambda}v(x_i)\chi_{\Omega_i},
$$
where $\chi_{\Omega_i}$ denotes the indicator function
of the set $\Omega_i.$

Due to stability reasons, especially for the case of dominating convection,
the class of finite volume methods under consideration
is characterized by an additional stabilization technique
called \emph{upwinding.}
For that purpose we introduce a scaling function
$K:\ \R\to [0,\infty)$ which is defined by the help
of a weighting function $r:\;\R\to [0,1]$ as $K(z):=1-[1-r(z)]z.$

A typical example of such  a weighting function is
\begin{equation}\label{r_ex2}
r(z):=1-\frac{1}{z}\Big(1-\frac{z}{e^z-1}\Big),
\end{equation}
leading to $K(z)=z/(e^z-1),$ the Bernoulli function.

The discrete problem for the case of a scalar diffusion coefficient,
i.e.\ where $\AAb$ is of the form $A\II$ with $A:\;\Omega\times\R\to\R$ and
$\II$ being the identity in $\R^d,$ is formulated as follows:
\begin{quote}
Find $\uT\in\VT$ such that
\end{quote}
\begin{equation}\label{eq12}
\forall \vT\in \VT:\quad a_\mathcal{T}(\uT;\vT)=\langle\fT,\vT\rangle,
\end{equation}
where
\begin{eqnarray*}
a_\mathcal{T}(\wT;\vT)&:=& \sumi \vTi \left\{\sumj
\dfrac{\mu_{ij}}{d_{ij}}K\left(\dfrac{\gamma_{ij}d_{ij}}%
{\mu_{ij}}\right)(\wTi-\wTj )m_{ij}+c_i \wTi m_i\right\},\\
\langle\fT,\vT\rangle&:=& \sumi f_i\vTi m_i,
%\quad\mbox{and}\quad
\end{eqnarray*}
and
\begin{eqnarray*}
\mu_{ij}=\mu_{ij}(\wTi,\wTj)
&:=&A\left(\frac{x_i+x_j}{2},\frac{\wTi+\wTj}{2}\right),\\
\gamma_{ij}=\gamma_{ij}(\wTi,\wTj)
&:=&\nu_{ij}\cdot \bb\left(\frac{x_i+x_j}{2},\frac{\wTi+\wTj}{2}\right),\\
c_i=c_i(\wTi)&:=& c(x_i,\wTi),
\quad
f_i:=f(x_i).
\end{eqnarray*}
Moreover, we introduce the following norms and seminorms
on $\VT:$
\begin{eqnarray*}
\|\vT\|_\mathcal{T}&:=& \sqrt{(\vT,\vT)_\mathcal{T}}
\ =\ \|\LT  \vT\|_{0,2,\Omega},\nonumber\\
|\vT|_V&:=& \left\{\sumi \vTi \sumj
(\vTi -\vTj )\dfrac{m_{ij}}{d_{ij}}\right\}^{1/2},
\label{V-seminorm}\\
\|\vT\|_V&:=& \left\{|\vT|_V^2+\|\vT\|_\mathcal{T}^2\right\}^{1/2}.\label{V-norm}
\end{eqnarray*}
For the sake of consistency in the notations, we also use
the following abbreviations of wellknown seminorms and norms in
the Sobolev space $H^1(\Omega)$:
$$
|\vT|_D:= |\vT|_{1,2,\Omega},\quad
\|\vT\|_D:= \|\vT\|_{1,2,\Omega}.
$$
The scheme (\ref{eq12}) with the weighting function (\ref{r_ex2})
is often called \emph{exponentially upwinded.}
It can be defined for other control functions $r:\R\to [0,1],$ too.
However, we have to assume that all of these control functions satisfy the
following properties:
\begin{description}
\item[\refsteppropl{rp:monot}\rm (\theproperty)]
$\quad r(z)$ is monotone for all $z\in\R,$ %old: P1
\item[\refsteppropl{rp:bnds}\rm (\theproperty)]
$\quad \lim\limits_{z\rightarrow -\infty}r(z)=0,\quad
\lim\limits_{z\rightarrow\infty}r(z)=1,$ %old: P2
\item[\refsteppropl{rp:mmatr}\rm (\theproperty)]
$\quad 1+zr(z)\ge 0\quad$ for all $z\in\R,$ %old: P3
\item[\refsteppropl{rp:sym}\rm (\theproperty)]
$\quad [1-r(z)-r(-z)]z=0\quad$ for all $z\in\R,$ %old: P4
\item[\refsteppropl{rp:possdef}\rm (\theproperty)]
$\quad \left[r(z)-\dfrac{1}{2}\right]z\ge 0\quad$ for all $z\in\R,$ %old: P5
\item[\refsteppropl{rp:cont}\rm (\theproperty)]
$\quad zr(z)$ is Lipschitz-continuous for all $z\in\R.$ %old: P6
\end{description}
We get from (\ref{rp:sym}) the relation
\begin{equation}\label{eqP3}
1+zr(z)=K(-z).
\end{equation}
Replacing in (\ref{eqP3}) the argument $z$ by $-z,$ (\ref{rp:mmatr})
immediately implies
\begin{description}
\item[\refsteppropl{rp:Kdef}\rm (\theproperty)] $\quad K(z)\ge 0\quad$
for all $z\in\R.$
\end{description}
\begin{example}\label{r_ex1}
The function
$$
r(z)=\frac{1}{2}[\mbox{sign}\,z +1],
$$
due to \cite{Baba:81}, has been investigated in
\cite{Risch:86}, \cite{Risch:90} for a linear equation in (\ref{genbvp}).
This scheme is called \emph{fully upwinded.}
\end{example}
The next two examples are simple \emph{approximations} of (\ref{r_ex2}).
\begin{example}\label{r_ex3}
$$
r(z)=\left\{\begin{array}{ccccr}
0&,&z&<& -m\\ \frac{z+m}{2m}&,&|z| &\le&  m\\ 1&,&z&>& m
\end{array}\right. ,\quad 0<m\le 8,
$$
\end{example}
\begin{example}\label{r_ex4}
$$
r(z)=\left\{\begin{array}{ccccr}
0&,&z&<& -m\\0.5&,&|z|&\le& m\\1&,&z&>& m
\end{array}\right. ,\quad 0\le m\le 2.
$$
This function violates property (\ref{rp:cont}).
\end{example}
\begin{example}\label{ex:r_Samarski}
The choice of the function
$$
r(z)=\frac{1}{2}\left[\frac{z}{2+|z|}+1\right]
$$
goes back to Samarskij \cite{Samarskij:65}.
\end{example}
\begin{example}\label{r_ex5}
In \cite{McCartin:83} the function
$$
r(z)=\frac{1}{2}[\tanh z+1]
$$
was proposed.
\end{example}
\begin{example}\label{r_ex6}
Finally, it can be taken the function
$$
r(z)=\left\{\begin{array}{ccc}
(1-\sigma )/2&,&z<0\\(1+\sigma )/2&,&z\geq 0
\end{array}\right.
\quad\mbox{with}\quad\sigma (z):=\max\left\{
0,1-\frac{2}{|z|}\right\},
$$
what corresponds to Ikeda's \emph{partial upwind scheme E}
\cite{Ikeda:83}.
\end{example}
Both for computational and theoretical reasons
it may be advisable, in a really nonlinear situation,
to choose differentiable control functions $r.$
In the sequel, if there is no special reference, we assume that
the scheme under consideration is defined for a general function $r$
that possesses the properties (\ref{rp:monot}) to (\ref{rp:cont}).

Finally we mention two equivalent representations of the form $a_\mathcal{T}.$
First we remember that the leading coefficient
$\dfrac{\mu_{ij}}{d_{ij}}K\left(\dfrac{\gamma_{ij}d_{ij}}{\mu_{ij}}\right)$
in $a_\mathcal{T}$ can be written, by the definition of $K,$ in the following manner:
$$
\frac{\mu_{ij}}{d_{ij}}K\left(\frac{\gamma_{ij}d_{ij}}{\mu_{ij}}\right)
=\frac{\mu_{ij}}{d_{ij}}\left\{1-\frac{\gamma_{ij}d_{ij}}{\mu_{ij}}
\left[1-r\left(\frac{\gamma_{ij}d_{ij}}{\mu_{ij}}\right)\right]\right\}
=\frac{\mu_{ij}}{d_{ij}}-\left(1-r_{ij}\right)\gamma_{ij}\,,
$$
where
$r_{ij}:=r\Big(\dfrac{\gamma_{ij}d_{ij}}{\mu_{ij}}\Big).$
Hence we get the representation
\begin{eqnarray}
&&a_\mathcal{T}(\wT;\vT)\nonumber\\
&=&\displaystyle\sumi \vTi \Big\{\sumj
\big\{\mu_{ij}(\wTi-\wTj )\dfrac{m_{ij}}{d_{ij}}
-\left(1-r_{ij}\right)(\wTi-\wTj )\gamma_{ij}m_{ij}\big\}\nonumber\\
&&\rule{14ex}{0ex} + c_i \wTi m_i\Big\}.\label{alrepr1}
\end{eqnarray}
Furthermore, introducing the notations
\begin{eqnarray}
a_\mathcal{T}^0(\wT;\vT)&:=&\sumi \vTi \sumj
\mu_{ij}(\wTi-\wTj )\dfrac{m_{ij}}{d_{ij}},\nonumber\\
b_\mathcal{T}(\wT;\vT)&:=&\sumi \vTi \sumj
\left[(1-r_{ij})\wTj -\left(\dfrac{1}{2}-r_{ij}\right)\wTi\right]
\gamma_{ij}m_{ij},\label{bl1}\\
d_\mathcal{T}(\wT;\vT)&:=&\sumi \left\{c_im_i
-\dfrac{1}{2}\sumj \gamma_{ij}m_{ij}\right\}\wTi \vTi ,
\label{bl2}
\end{eqnarray}
we get a splitting of $a_\mathcal{T}$ which is comparable with (\ref{asplit}):
\begin{equation}\label{a_lsplit}
a_\mathcal{T}(\wT;\vT)=a_\mathcal{T}^0(\wT;\vT)+b_\mathcal{T}(\wT;\vT)+d_\mathcal{T}(\wT;\vT).
\end{equation}
\begin{remark}\label{divappr_rem}
In the special case $\nabla\cdot \bb\equiv 0$ on $\Omega,$
it is senseful to use the following versions of $b_\mathcal{T}$
and $d_\mathcal{T}$:
\begin{eqnarray*}
b_\mathcal{T}(\wT;\vT)&=&\sumi \vTi \sumj
\left[(1-r_{ij})\wTj +r_{ij}\wTi\right]\gamma_{ij}m_{ij},\\
d_\mathcal{T}(\wT;\vT)&=&\sumi c_i\wTi \vTi m_i.
\end{eqnarray*}
\end{remark}
\subsection{The case of Voronoi diagrams and matrix-valued diffusion coefficients}
Using the representation
$$
\wT=\sum_{j\in\Lambda_T} \wTj \psi_j
$$
on a single element $T,$ where $\{\psi_j\}_{j\in\Lambda}$
is the standard nodal basis of $\VT,$
we easily see that
\begin{eqnarray*}
\int_{\partial\Omega_i\cap T}(\AAb\nabla \wT)\cdot\nu\,ds
&=&\sum_{j\in\Lambda_T}
\int_{\partial\Omega_i\cap T}\wTj (\AAb\nabla \psi_j)\cdot\nu\,ds\\
&=&\sum_{j\in\Lambda_T\setminus\{i\}} \left( \wTj  - \wTi \right)
\int_{\partial\Omega_i\cap T}(\AAb\nabla \psi_j)\cdot\nu\,ds\,.
\end{eqnarray*}
In the next step the matrix $\AAb$ is approximated
by a piecewise constant matrix $\AAb_\mathcal{T}.$
Summing up over all elements $T$ which lie in the support of $\psi_i,$
we have the relation
$$
\int_\Omega(\AAb_\mathcal{T} \nabla \wT)\cdot\nabla\psi_i\,dx
=\sum_{j\in\Lambda_i} \left( \wTi - \wTj  \right)
\int_{\partial\Omega_i}(\AAb_\mathcal{T} \nabla \psi_j)\cdot\nu\,ds\,.
$$
With the definition
\begin{equation}\label{eq:mvmudef}
\mu_{ij}:=
\left\{\begin{array}{l@{\quad}l}
\displaystyle\frac{d_{ij}}{m_{ij}}
\int_{\partial\Omega_i}(\AAb\nabla \psi_j)\cdot\nu\,ds,&m_{ij}>0\,,\\[.5ex]
0\,,&m_{ij}=0\,,
\end{array}\right.
\end{equation}
it follows that
$$
\int_\Omega(\AAb\nabla \wT)\cdot\nabla\psi_i\,dx
\approx\int_\Omega(\AAb_\mathcal{T}\nabla \wT)\cdot\nabla\psi_i\,dx
=\sum_{j\in\Lambda_i}\mu_{ij}\left( \wTi - \wTj  \right)
\frac{m_{ij}}{d_{ij}}\,.
$$
Unfortunately, it is wellknown that in the case $d=3$ the right equality does not hold.
Nevertheless, the right-hand side -- together with the above definition (\ref{eq:mvmudef})
of $\mu_{ij}$ -- is senseful for $d=3,$ and thus this formula can be used for discretization.

Consequently, in order to obtain a discretization for the case
of a matrix-valued diffusion coefficient, it is sufficient
to replace in the forms $a_\mathcal{T}^0$ and $b_\mathcal{T}$
the corresponding values of $\mu_{ij}$ according to formula (\ref{eq:mvmudef}).

\begin{remark}
The really critical point in the discretization of diffusion-convec\-tion equations
with matrix-valued diffusion coefficients
consists in the appropriate choice of the stabilization mechanism
in the situation where the eigenvalues of $\AAb$ are widely spreaded
(cf. \cite{Angermann:00c}, \cite{Angermann:05b}).
\end{remark}
\subsection{The case of Donald diagrams}
Let us now consider a family of admissible (in the sense of FEM,
cf.\ \cite[Ch.\ 2]{Ciarlet:78}) triangulations
$\mathcal{F}=\{\mathcal{T}\}.$
Then, for any $T\in\mathcal{T}$ with local vertices
$z_j\equiv x_{i_j},\ i_j\in\Lambda_T,\ j\in[1,d+1]_{\N},$
we define
$$
\Omega_{i_j,T}^D:=\left\{x\in T:\left(\forall k\in [1,d+1]_{\N}\setminus\{j\}:
\lambda_k(x)<\lambda_j(x)\right)\right\},
$$
where $\lambda_j(x)$ is the $j$-th barycentric coordinate of $x$ w.r.t.\ $T.$
Define for $i\in\overline\Lambda$ the sets
$$
\Omega_i^D:=\inte\left(\bigcup_{T:\,\partial T\ni x_i}
\overline{\Omega_{i,T}^D}\right).
$$
In this way, we get a family of Donald diagrams.

Although it is possible to introduce a discretization
like (\ref{a_lsplit}), we use the following version:
\begin{equation}\label{barydiscr}
a_\mathcal{T}(\wT;\vT)=(\AAb(\cdot,\wT)\nabla \wT,\nabla \vT)
+b_\mathcal{T}(\wT;\vT)+d_\mathcal{T}(\wT;\vT),
\end{equation}
where the forms $b_\mathcal{T},\,d_\mathcal{T}$ are defined analogously to
(\ref{bl1}),(\ref{bl2}).
In particular, $\gamma_{ij}\in\R$ is an approximation to
$(\nu\cdot \bb)(\cdot,\wT)|_{\Gamma{ij}}.$

In the case of a matrix-valued diffusion coefficient, we define $\mu_{ij}$ analogously
to (\ref{eq:mvmudef}) but use it only in $b_\mathcal{T}$ to ensure
a certain stabilization. The form $a_\mathcal{T}^0$ remains as it is,
i.e. $a_\mathcal{T}^0(\wT;\vT):=(\AAb(\cdot,\wT)\nabla \wT,\nabla \vT).$
\section{Stability and a priori error estimates}
\label{sec:aprioriresults}
\subsection{The case of a linear equation with a scalar diffusion coefficient}
In this section we give a short review of some wellknown
properties of the schemes (\ref{eq12}) and (\ref{barydiscr})
for the case of a linear equation with a scalar diffusion coefficient.
We start with the formulation of conditions
with respect to the approximations $\mu_{ij}$ and
$\gamma_{ij}.$
\begin{description}
\item[\rm (A2.1)] $\mu_{ij}$ is an approximation of
the term $m_{ij}^{-1}\int_{\Gamma_{ij}}A\,ds$ satisfying the
following conditions:
\begin{enumerate}
\renewcommand{\labelenumi}{(\roman{enumi})}
\item $0\le\mu_{ij}\le \|A\|_{1,\infty,\Omega},$
\item $\mu_{ij}=\mu_{ji},$
\item $\left|\mu_{ij}-m_{ij}^{-1}\int_{\Gamma_{ij}}A\,ds\right|
\le Ch_T|A|_{1,\infty,\Omega},$
where $T$ is one of the simplices having the vertices $x_i,x_j,$
and $C>0$ is a constant independent of $a,h_T,i,j.$
\end{enumerate}
\item[\rm (A2.2)] $\gamma_{ij}$ is an approximation of
the term $m_{ij}^{-1}\int_{\Gamma_{ij}}\nu\cdot \bb\,ds$ satisfying the
following conditions:
\begin{enumerate}
\renewcommand{\labelenumi}{(\roman{enumi})}
\item $|\gamma_{ij}|\le \|\bb\|_{1,\infty,\Omega},$
\item $\gamma_{ij}=-\gamma_{ji},$
\item $\left|\gamma_{ij}
-m_{ij}^{-1}\int_{\Gamma_{ij}}(\nu_{ij}\cdot \bb)ds\right|
\le Ch_T|\bb|_{1,\infty,\Omega},$
where $T$ is one of the simplices having the vertices $x_i,x_j,$
and $C>0$ is a constant independent of $\bb,h_T,i,j.$
\end{enumerate}
\end{description}
The subsequent results are extensions of the theory developed in
\cite{Angermann:91a}, \cite{Angermann:95b},
see also \cite[Ch.\ 6]{Angermann:03d}.
\begin{theorem}[Discrete coercivity]\label{coral-coerc}
Let a family $\mathcal{F}=\{\mathcal{T}\}$ of triangulations
be given, where in the special case of Voronoi diagrams
(i.e.\ $\Xi=V$) all elements $T$ are self-centered
and in the special case of Donald diagrams
(i.e.\ $\Xi=D$) the family is shape-regular.
Moreover, let the assumptions (A2.1), (A2.2) be satisfied.
Then, for $h_0>0$ sufficiently small there exist two constants
$\overline{a}_0>0$ and $\overline{a}_1>0$ independent of $h$ such
that for all $h\in (0,h_0]$ and $\vT\in \VT$ the relation
$$
a_\mathcal{T}(\vT;\vT)\ge\overline{a}_0|\vT|_\Xi^2
+\overline{a}_1\|\vT\|_\mathcal{T}^2
$$
holds.
\end{theorem}

The a priori error estimate is based on this stability property
and on the following consistency result.
\begin{lemma}[Discrete consistency]\label{weakconsist}
Let a shape-regular family $\mathcal{F}$ of triangulations $\{\mathcal{T}\}$
be given, where in the special case of Voronoi diagrams
(i.e.\ $\Xi=V$) all elements $T$ are self-centered,
and let the assumptions (A2.1), (A2.2) be satisfied.
Then, if $h_0>0$ is sufficiently small, for any element
$w\in W_2^2(\Omega)\cap V$ and any element $\vT\in \VT$
the estimate
\begin{eqnarray*}
&&\left|a_\mathcal{T}(\IT w;\vT)
-(-\nabla\cdot(A\nabla w)+\bb\cdot\nabla w +cw,\LT  \vT)\right|\\
&\le& Ch\|w\|_{2,2,\Omega}\left[|\vT|_\Xi+\|\vT\|_\mathcal{T}\right]
\end{eqnarray*}
holds for all $h\in (0,h_0],$ where $C>0$ is a constant which does not
depend on $h.$
\end{lemma}
The proof of the following theorem is a modification of the standard
proof of Strang's first lemma.
\begin{theorem}[A priori error estimate]\label{errortheorem}
Let a shape-regular family $\mathcal{F}$ of triangulations $\{\mathcal{T}\}$
be given, where in the special case of Voronoi diagrams
(i.e.\ $\Xi=V$) all elements $T$ are self-centered,
let the assumptions
(A2.1), (A2.2)
be satisfied and suppose that the solution $u\in V$ of problem (\ref{genbvp})
additionally belongs to $W_2^2(\Omega).$

Then, for sufficiently small $h_0>0$ the estimate
$$
\|u-\uT\|_\Xi
\le Ch\left[\|u\|_{2,2,\Omega}+|f|_{1,q,\Omega}\right]
$$
holds for all $h\in (0,h_0]$, where the constant $C>0$ is independent of $h.$
\end{theorem}

\subsection{The quasilinear case}

Due to the  possible structural diversity of the nonlinearities
in (\ref{genbvp}),
in the nonlinear situation there is not such a
relatively canonical theory as in the linear case.

We mention here only a few papers which are concerned with
the investigation of node-centered finite volume methods for nonlinear
elliptic (or parabolic) equations and refer to the literature
cited therein:
\cite{Fuhrmann:01a},
\cite{Chatzipantelidis:05b},
\cite{Eymard:06}.

\section{\emph{A posteriori} error estimates for nonlinear problems}
In this section we present the general approach that does not depend
on the particular discretization.

The nonlinear primal problem we are interested in is given by
\begin{equation}\label{eq:nonlinproblem}
u \in V:\ a(u;v) + a_\delta(u;v) = \langle f,v \rangle \quad\forall v\in V.
\end{equation}
It represents the weak formulation of the originally given (accurate)
boundary-value problem
for a partial differential equation in a real Hilbert space $V,$
where $f$ is a linear functional on $V$ and $\langle f,v \rangle$
denotes the value of $f$ at the element $v\in V.$
The forms $a:\;V\times V\to\R$ and $a_\delta:\;V\times V\to\R$
are linear in the second argument but may be nonlinear in the first one.
In the context of the boundary-value problem (\ref{genvar}),
the left-hand side of (\ref{genvar}) is written in (\ref{eq:nonlinproblem})
as the sum $a+a_\delta,$ where $a$ stands for a certain simplified problem
and $a_\delta$ represents a part of the equation
which is to be neglected in the practical computations.
That is, the discretization applies only to $a$ in (\ref{eq:nonlinproblem}).
The goal is to estimate the influence of both neglecting $a_\delta$
and discretizing $a$ and $f$ with respect to a given output functional
$j:\;V\to\R.$
\begin{example}
Consider (\ref{genbvp}) with
$$
\AAb(x,w):=\eps(x)|w|^{\gamma(x)}\II,
\quad
\bb(x,w):=\bb_0(x)|w|^{\gamma(x)/2},
\quad
c(x,w):=c_0(x),
$$
where $\eps, c_0, f, \gamma:\;\Omega\to\R$ and $\bb_0:\;\Omega\to\R^d$
are smooth functions (satisfying certain additional conditions,
in particular $-1<\gamma_-\le\gamma(x)\le\gamma_+<\infty$ on $\Omega$
for some constants $\gamma_-,\gamma_+$).
Then, for some constant elements $w_0, \gamma_0\in\R,$ we can set
$$
a(w;v):=(\eps|w_0|^{\gamma_0}\nabla w,\nabla v)
+(|w_0|^{\gamma_0/2}\bb_0\cdot\nabla w,v)
+(c_0w,v),
$$
and $a_\delta(w;v)$ is the canonical error term with respect to the
correct weak formulation of (\ref{genbvp}).
\end{example}
The directional derivatives of $a(u;\cdot)$ and $a_\delta(u;\cdot)$
in $u$ will be denoted by $a'(u;\cdot,\cdot)$
and
$a_\delta'(u;\cdot,\cdot),$ respectively. The form
$$
a'(u;w,v) := \lim_{\eps\to 0}\frac{1}{\eps}\left[a(u + \eps w;v) - a(u;v)\right]
$$
is linear in $w$ and $v.$
The second and third directional derivatives are denoted by
$a''(u;\cdot, \cdot, \cdot)$
and
$a'''(u;\cdot, \cdot, \cdot, \cdot),$ respectively.
In the general case of a nonlinear output functional $j,$
the corresponding dual problem we will use in the analysis is the following:
\begin{equation}\label{eq:nonlindualproblem}
z \in V:\ a'(u;w,z) + a_\delta'(u;w,z) = j'(u;w) \quad\forall w\in V.
\end{equation}
The solution $z\in V$ of the dual problem is called \emph{influence function}
for the particular choice of $j$ (\cite{Ainsworth:00}).
The primal solution $u_m \in V$ and the dual solution $z_m\in V$
of the reduced problems are given by
\begin{eqnarray}
u_m \in V :\ a(u_m;v) &=& \langle f,v \rangle \quad\forall v\in V,
\label{eq:rednonlinproblem}\\
z_m \in V :\ a'(u_m;w, z_m) &=& j'(u_m;w) \quad\forall w\in V.
\label{eq:rednonlindualproblem}
\end{eqnarray}
These variational problems will be formulated in terms of optimization problems.
The primal and dual solutions will be expressed by the variables
$x := (u, z) \in X := V\times V$ and $x_m := (u_m, z_m) \in X.$
In the variational space $X,$ we consider the functionals
\begin{eqnarray}
L(x) &:=& L_m(x) + L_\delta(x),\label{eq:defL}\\
L_m(x) &:=& j(u) + \langle f,z \rangle - a(u;z),\label{eq:defLm}\\
L_\delta(x) &:=& -a_\delta(u;z).\label{eq:defdeltaL}
\end{eqnarray}
The derivative of $L$ applied to a test function $y = (w,v)\in X$ is
$$
L'(x;y) = j'(u;w) - a'(u;w, z) - a_\delta'(u;w,z)
+ \langle f, v\rangle - a(u;v) - a_\delta(u;v).
$$
Obviously, the original primal and dual problems
(\ref{eq:nonlinproblem}) and (\ref{eq:nonlindualproblem})
and the reduced primal and dual problems
(\ref{eq:rednonlinproblem}) and (\ref{eq:rednonlindualproblem})
consist of finding the stationary points
$x = (u, z)$ and $x_m = (u_m, z_m)$ of $L$ and $L_m,$ respectively:
\begin{eqnarray}
x \in X:\ L'(x;y) &=& 0 \quad\forall y \in X,\label{eq:statL}\\
x_m \in X:\ L_m'(x_m;y) &=& 0 \quad\forall y \in X.\label{eq:statLm}
\end{eqnarray}
Furthermore, the target quantities are given by evaluation of $L$ and $L_m$
at the following stationary points:
$$
j(u) = L(x),
\quad
j(u_m) = L_m(x_m).
$$
In order to balance the model and discretization errors,
we have to include the discretization error in the %a posteriori
analysis.
To do this, let $\VT\subset V$ be a finite-dimensional subspace.
Typically $\VT$ is a finite element space with respect to
a partition $\mathcal{T}$ of the computational domain
$\Omega\subset\R^d,$ $d\in\{2,3\},$ where possible homogeneous Dirichlet
boundary conditions are already included in
the choice of the spaces $V$ and $\VT.$
Let $a_\mathcal{T}:\;\VT\times\VT\to\R$ be a nonlinear form
which is different, in general, from the simple restriction of $a$
to $\VT\times\VT,$
and denote by $f_\mathcal{T}:\;\VT\to\R$ a linear functional
which not necessarily coincides with $f|_{\VT}.$
For instance, $a_\mathcal{T}$ and $f_\mathcal{T}$ may result
from the finite volume discretization of $a,$ $f$ in (\ref{eq:nonlinproblem})
according to Section \ref{s:fv_deriv}.

Then $u_{\mathcal{T}m}\in\VT$ is the discrete solution of the problem
\begin{equation}\label{eq:fvmrednonlindualproblem}
u_{\mathcal{T}m} \in \VT:\ a_\mathcal{T}(u_{\mathcal{T}m};v)
= \langle f_\mathcal{T},v \rangle
\quad\forall v\in \VT
\end{equation}
involving both types of error.
The operators $L$ and $L_m$ are still given by (\ref{eq:defL})-(\ref{eq:defdeltaL}).
The difference lies in the definition of the discrete solution
$x_{\mathcal{T}m} = (u_{\mathcal{T}m}, z_{\mathcal{T}m})
\in X_\mathcal{T} := \VT\times \VT,$
where now $u_{\mathcal{T}m}$ satisfies (\ref{eq:fvmrednonlindualproblem})
and $z_{\mathcal{T}m}$ is the solution of the following dual problem:
\begin{equation}\label{eq:discrrednonlindualproblem}
z_{\mathcal{T}m} \in \VT:\ a'(u_{\mathcal{T}m};w, z_{\mathcal{T}m})
= j'(u_{\mathcal{T}m};w) \quad\forall w\in \VT.
\end{equation}
In such a setting, the relations
$\ a(u_{\mathcal{T}m};v) = \langle f,v \rangle \ $
and $\ L_m'(x_{\mathcal{T}m};y) = 0 \ $
are no longer valid for all $v\in \VT$ resp. $y \in X_\mathcal{T}.$

The target quantities are given by the evaluation of $L$ and $L_{\mathcal{T}m},$
where
\begin{equation}\label{eq:defdeltaLfv}
L_{\mathcal{T}m}(x):=j(u) + \langle f_\mathcal{T},z \rangle - a_\mathcal{T}(u;z),
\end{equation}
at the following stationary points:
\begin{equation}\label{eq:fvmtargets}
j(u) = L(x),
\quad
j(u_m) = L_{\mathcal{T}m}(x_m).
\end{equation}
For the formulation of the error representation, we use the following notation
for the primal and dual residual with respect to the reduced model and for test functions
$(w,v) \in X$:
\begin{eqnarray*}
\varrho(u_{\mathcal{T}m};v) &:=& \langle f,v \rangle - a(u_{\mathcal{T}m};v),\\
\varrho^*(u_{\mathcal{T}m};z_{\mathcal{T}m},w)
&:=& j'(u_{\mathcal{T}m};w) - a'(u_{\mathcal{T}m};w,z_{\mathcal{T}m}).
\end{eqnarray*}
\begin{theorem}\label{th:discrerr}
If $a(u;\cdot),$ $a_\delta(u;\cdot)$ and the functional $j(u)$ are sufficiently
differentiable with respect to $u,$ then we have
\begin{eqnarray*}
j(u) - j(u_{\mathcal{T}m}) &=& -a_\delta(u_{\mathcal{T}m};z_{\mathcal{T}m})\\
&& +\  \langle f,z_{\mathcal{T}m} \rangle - \langle f_\mathcal{T},z_{\mathcal{T}m} \rangle
- a(u_{\mathcal{T}m};z_{\mathcal{T}m}) + a_\mathcal{T}(u_{\mathcal{T}m};z_{\mathcal{T}m})\\
&&+\ \frac{1}{2}\left[\varrho(u_{\mathcal{T}m};z - i_\mathcal{T} z)
+ \varrho^*(u_{\mathcal{T}m};z_{\mathcal{T}m},u - i_\mathcal{T} u)\right]\\
&&-\ \frac{1}{2} \left[a_\delta(u_{\mathcal{T}m};e_z)
+ a_\delta'(u_{\mathcal{T}m};e_u, z_{\mathcal{T}m})\right]\\
&&-\ \frac{1}{2}\,\varrho(u_{\mathcal{T}m};z_{\mathcal{T}m}- i_\mathcal{T} z)
-\frac{1}{2}\,R,
\end{eqnarray*}
where $e := (e_u, e_z) := (u - u_{\mathcal{T}m}, z - z_{\mathcal{T}m}),$
$i_\mathcal{T}:\ V\to \VT$ is an interpolation operator,
and the remainder $R$ is given by
$$
R := \int_0^1 \sigma(1-\sigma)L'''(x_{\mathcal{T}m} + \sigma e;e, e, e)\,d\sigma.
$$
\end{theorem}
\proof
By (\ref{eq:fvmtargets}),
\begin{eqnarray*}
j(u) - j(u_{\mathcal{T}m})
&=& L(x) - L_{\mathcal{T}m}(x_{\mathcal{T}m})\\
&=& L(x) - L_m(x_{\mathcal{T}m}) + L_m(x_{\mathcal{T}m}) - L_{\mathcal{T}m}(x_{\mathcal{T}m})\\
&=& L(x) - L_m(x_{\mathcal{T}m})\\
&& +\  \langle f,z_{\mathcal{T}m} \rangle - \langle f_\mathcal{T},z_{\mathcal{T}m} \rangle
- a(u_{\mathcal{T}m};z_{\mathcal{T}m}) + a_\mathcal{T}(u_{\mathcal{T}m};z_{\mathcal{T}m}),
\end{eqnarray*}
where the last step is a consequence of the definitions (\ref{eq:defdeltaL}), (\ref{eq:defdeltaLfv}).

The first difference can be estimated as in the proof of \cite[Thm.\ 2.1]{Braack:03}:
\begin{eqnarray*}
L(x) - L_m(x_{\mathcal{T}m}) &=& L(x) - L(x_{\mathcal{T}m}) + L_\delta(x_{\mathcal{T}m})\\
&=& \int_0^1 L'(x_{\mathcal{T}m} + \sigma(x-x_{\mathcal{T}m});x-x_{\mathcal{T}m})d\sigma
 + L_\delta(x_{\mathcal{T}m})\\
&=& \frac{1}{2}\left[L'(x_{\mathcal{T}m};e)+L'(x;e) - R\right]
- a_\delta(u_{\mathcal{T}m};z_{\mathcal{T}m})
\end{eqnarray*}
with the above given remainder $R$ of the trapezoidal rule.
Since $L'(x;e)=0$ by (\ref{eq:statL}), we get
\begin{eqnarray*}
L(x) - L_m(x_{\mathcal{T}m}) &=& -a_\delta(u_{\mathcal{T}m};z_{\mathcal{T}m})
+ \frac{1}{2} \left[L'(x_{\mathcal{T}m};e) - R\right].
\end{eqnarray*}
Furthermore,
\begin{eqnarray*}
L'(x_{\mathcal{T}m};e) &=& %L_m'(x_{\mathcal{T}m};e) + L_\delta'(x_{\mathcal{T}m};e)\\
j'(u_{\mathcal{T}m};e_u) - a'(u_{\mathcal{T}m};e_u,z_{\mathcal{T}m})
- a_\delta'(u_{\mathcal{T}m};e_u,z_{\mathcal{T}m})\\
&&+\ \langle f, e_z\rangle - a(u_{\mathcal{T}m};e_z) - a_\delta(u_{\mathcal{T}m};e_z)\\
&=& %\varrho(u_{\mathcal{T}m};e_z)
\varrho^*(u_{\mathcal{T}m};z_{\mathcal{T}m},e_u)
- a_\delta'(u_{\mathcal{T}m};e_u, z_{\mathcal{T}m})
+\varrho(u_{\mathcal{T}m};e_z) - a_\delta(u_{\mathcal{T}m};e_z).
\end{eqnarray*}
Since the Galerkin orthogonality is violated, in general,
we cannot use the standard argument
$$
0 = \varrho(u_{\mathcal{T}m};z_{\mathcal{T}m})
= \varrho(u_{\mathcal{T}m};i_\mathcal{T} z)
$$
to replace $z_{\mathcal{T}m}$ by $i_\mathcal{T} z$ in the third term.
Here we can only make use of an analogous property of the dual problem
(\ref{eq:discrrednonlindualproblem}), i.e.
$$
0 = \varrho^*(u_{\mathcal{T}m};z_{\mathcal{T}m},u_{\mathcal{T}m})
= \varrho^*(u_{\mathcal{T}m};z_{\mathcal{T}m},i_\mathcal{T} u).
$$
(Of course, if the dual problem is approximated by a finite volume method, too,
then we have to argue as for the primal problem.)
Thus we arrive at
\begin{eqnarray}
L'(x_{\mathcal{T}m};e)
&=& \varrho^*(u_{\mathcal{T}m};z_{\mathcal{T}m},u - i_\mathcal{T} u)
+\varrho(u_{\mathcal{T}m};z- i_\mathcal{T} z)\label{eq:ezsplit}\\
&&-\ \varrho(u_{\mathcal{T}m};z_{\mathcal{T}m}- i_\mathcal{T} z)
- a_\delta(u_{\mathcal{T}m};e_z) - a_\delta'(u_{\mathcal{T}m};e_u, z_{\mathcal{T}m}).
\nonumber
\end{eqnarray}
This gives the assertion.
\close

In order to use numerically the error representation derived in Theorem \ref{th:discrerr},
we have to approximate various terms.
In particular, we will neglect the higher-order terms in $e,$
namely the remainder $R$ and the terms $a_\delta(u_{\mathcal{T}m};e_z),$
$a_\delta'(u_{\mathcal{T}m};e_u, z_{\mathcal{T}m}),$
cf.\ the related discussion in \cite{Braack:03}.
%with respect to the model $a_\delta(\cdot;\cdot).$
Furthermore, we have to approximate the interpolation errors
$u - i_\mathcal{T}u$ and $z - i_\mathcal{T}z.$
An efficient possibility for doing this is the recovery process
of the computed quantities by patch-wise higher-order interpolation
expressed via the operator
$i_\mathcal{T}^+:\;\VT\to\VT^+$ formally,
where $\VT^+$ is a richer discrete space than $\VT$
(see \cite[Sect.\ 5]{Becker:01b}, \cite[Sect.\ 3.2]{Rannacher:05}).
For instance, in the case of triangles ($d = 2$) or tetrahedra ($d = 3$)
and when $\VT$ consists of piecewise linear elements,
quadratic interpolation may be used.
For quadrilaterals and piecewise $d$-linear elements,
the interpolation can be done on $d$-quadratic elements.
In order to preserve a sufficient high accuracy of the interpolation
procedure, special care on elements with hanging nodes is required.

The interpolation errors will be numerically approximated by
\begin{eqnarray}
z - i_\mathcal{T} z &\approx&  i_\mathcal{T}^+ z_{\mathcal{T}m}
- z_{\mathcal{T}m},\label{eq:ezapprox}\\
u - i_\mathcal{T} u &\approx&  i_\mathcal{T}^+ u_{\mathcal{T}m}
- u_{\mathcal{T}m}.\nonumber
\end{eqnarray}
Without the modeling error and in the case of conforming methods,
this approximation is usually observed to be accurate enough.

Taking into account that the residual
$\varrho^*(u_{\mathcal{T}m};z_{\mathcal{T}m},v)$ vanishes
with respect to a discrete test function $v\in \VT,$
we obtain from Theorem \ref{th:discrerr}
the following approximate estimator consisting of three indicators:
\begin{eqnarray}
j(u) - j(u_{\mathcal{T}m}) &\approx& \eta_\mathcal{T} + \eta_m + \eta_{nc},\nonumber\\
\eta_\mathcal{T}
&:=& \frac{1}{2} \left[
\varrho(u_{\mathcal{T}m};i_\mathcal{T}^+ z_{\mathcal{T}m}- z_{\mathcal{T}m})
+ \varrho^*(u_{\mathcal{T}m};z_{\mathcal{T}m},i_\mathcal{T}^+ u_{\mathcal{T}m})\right],
\label{eq:partsest}\\
\eta_m &:=& -a_\delta(u_{\mathcal{T}m};z_{\mathcal{T}m}),\nonumber\\
\eta_{nc}&:=& \langle f,z_{\mathcal{T}m} \rangle
- \langle f_\mathcal{T},z_{\mathcal{T}m} \rangle
- a(u_{\mathcal{T}m};z_{\mathcal{T}m})
+ a_\mathcal{T}(u_{\mathcal{T}m};z_{\mathcal{T}m}).\nonumber
\end{eqnarray}
The indicator $\eta_\mathcal{T}$ of the approximate estimator can be considered as
the conforming contribution of the discretization,
and the indicator $\eta_m$ measures the influence of the model.
For complex models, the evaluation of $\eta_m$ may be expensive.
Often in practice the decomposition $a+a_\delta$ is changed successively
in such a way that portions of $a_\delta$ are (locally) shifted to $a.$
The indicator $\eta_{nc}$ results from the nonconformity of the discretization
method caused by the violation of the Galerkin orthogonality.
The practical treatment of  $\eta_{nc}$ will be discussed
in Section \ref{s:fvmest}.
\begin{remark}
According to (\ref{eq:ezapprox}), there are two ways for the treatment
of the term
$$
\varrho(u_{\mathcal{T}m};z - i_\mathcal{T} z)
-\varrho(u_{\mathcal{T}m};z_{\mathcal{T}m}- i_\mathcal{T} z)
$$
occuring in Theorem \ref{th:discrerr}.
Either we write it as $\varrho(u_{\mathcal{T}m};z - z_{\mathcal{T}m})$
(i.e.\ we reverse the splitting used in (\ref{eq:ezsplit}))
and replace then $z$ by $i_\mathcal{T}^+ z_{\mathcal{T}m},$
or we replace $z - i_\mathcal{T} z$ by
$i_\mathcal{T}^+ z_{\mathcal{T}m}- z_{\mathcal{T}m}$
and $i_\mathcal{T} z$ by $z_{\mathcal{T}m}.$
In both cases, we arrive at the same result:
$$
\varrho(u_{\mathcal{T}m};z - i_\mathcal{T} z)
-\varrho(u_{\mathcal{T}m};z_{\mathcal{T}m}- i_\mathcal{T} z)
\approx  \varrho(u_{\mathcal{T}m};i_\mathcal{T}^+ z_{\mathcal{T}m} - z_{\mathcal{T}m}).
$$
\end{remark}
In order to use the information (\ref{eq:partsest}) for changing locally the
model or the discretization parameters (e.g.\ the mesh size),
we have to localize the indicators.
After that, an adaptive process has to be designed
in order to balance the error sources.

Regarding the localization of $\eta_\mathcal{T}$ and $\eta_m,$
so here there are no new aspects. We refer, for instance, to \cite{Braack:03}.

\section{Application to the finite volume method}
\label{s:fvmest}
In the papers \cite{Angermann:91a}, \cite{Angermann:92a},
an extension of Babu\v{s}ka\&Rhein\-boldt's \emph{a posteriori} error estimates
for finite element methods (\cite{Babuska:78}) to finite volume methods
for linear diffusion-convection eqations
has been proposed.
In a subsequent paper (\cite{Angermann:95d}), for a singularly perturbed model
problem a modification was introduced with the aim to get two-sided bounds
of the error such that the constants occuring in these bounds are independent
of the perturbation parameter.
In \cite{Angermann:98a} and \cite{Thiele:99}, residual type error estimates
for finite volume discretizations of more complicated problems
in two and three space dimensions have been presented.
A rather general framework for the a posteriori estimation
in various finite volume methods can be found in \cite{Vohralik:08},
however this paper is restricted to linear problems and estimates
w.r.t.\ the energy norm.
In \cite{Angermann:10b}, \emph{dual-weighted residual error estimators}
for finite volume discretizations of linear diffusion-convection eqations
have been described. Here we apply the results of the previous section
to the nonlinear diffusion-convection problem.
As a result, we get \emph{a posteriori} estimates
for errors of functionals depending nonlinearly on the solution
and for possible modeling errors.

Interpreting $a_\mathcal{T}$ and $\fT$ as the finite volume discretizations (\ref{eq12})
of the forms $a$ and $f$ in (\ref{eq:nonlinproblem}),
we first observe that the the estimators $\eta_\mathcal{T}$ and $\eta_m$
depend only on the computed discrete solution but not directly on the structure
of $a_\mathcal{T}$ and $\fT.$
Therefore, these estimators can be treated as in the (conforming) finite element case
and we concentrate on the estimator $\eta_{nc}.$
To simplify the presentation, we will write $\xT=(\uT,\zT)$
instead of
$x_{\mathcal{T}m} = (u_{\mathcal{T}m}, z_{\mathcal{T}m}).$

Then, by definition, we have that
\begin{eqnarray}
\langle f,\zT \rangle-\langle f_\mathcal{T},\zT \rangle
&=&\sum_{T\in\mathcal{T}}\left\{(f,\zT)_T-(f,\zT)_{l,T}\right\}\label{eq:fldeco}\\
&:=&\sum_{T\in\mathcal{T}}\Big\{(f,\zT)_T-\sum_{i\in\Lambda_T}f_i\zTi m_i^T\Big\},
\nonumber
\end{eqnarray}
where $\displaystyle (f,\zT)_T:=\int_T f\zT\,dx$ and $m_i^T:=\measd{\Omega_i\cap T}.$
Analogously, with
\begin{eqnarray*}
&&a_{\mathcal{T},T}(\uT;\zT) \\
&:=&\sum_{i \in \Lambda} \zTi \left\{
\sum_{j\in\Lambda_T\setminus\{i\}} \!\!
\left\{ \mu_{ij} \frac{\uTi - \uTj }{d_{ij}}
-\gamma_{ij} \left(1-r_{ij}\right)(\uTi-\uTj ) \right\}m_{ij}^T
+ c_i \uTi m_i^T \! \right\}
\end{eqnarray*}
and $a_T(\uT;\zT)$ resulting from the restriction of all integrals
occuring in the expression for $a(\uT;\zT)$ to the domain of integration $T,$
we have that
\begin{equation}\label{eq:aldeco}
a_\mathcal{T}(\uT;\zT)-a(\uT;\zT)
=\sum_{T\in\mathcal{T}}\left\{a_{\mathcal{T},T}(\uT;\zT)-a_T(\uT;\zT)\right\}.
\end{equation}
Putting (\ref{eq:fldeco}) and (\ref{eq:aldeco}) together,
we conclude that
$$
\eta_{nc}
= \sum_{T\in\mathcal{T}}\Big\{(f,\zT)_T-\sum_{i\in\Lambda_T}f_i\zTi m_i^T\Big\}
+\sum_{T\in\mathcal{T}}\left\{a_{\mathcal{T},T}(\uT;\zT)-a_T(\uT;\zT)\right\}.
$$
This is the starting point for the practical, localized computation.

\section{Analysis of the nonconformity estimator}
In this section we show for the case of a scalar diffusion
coefficient $\AAb=A\II$ with $A:\;\Omega\times\R\to\R$
that $\eta_{nc}$ is \emph{order-consistent}
with the \emph{a priori} error estimate (Theorem \ref{errortheorem}).
A precise formulation of this property is given at the end of the section.

Using (\ref{eq:fldeco}), (\ref{eq:aldeco}), we get
the following decomposition:
\begin{eqnarray*}
\eta_{nc}
&=&\sum_{i\in\overline\Lambda} \int_{\Omega_i}[f\zT-f_i\zTi ]dx\\
&&+\sumi \zTi \sumj\mu_{ij}(\uTi-\uTj )\dfrac{m_{ij}}{d_{ij}}-(A\nabla \uT, \nabla\zT)\\
&&+\sumi \sumj (1-r_{ij})\gamma_{ij}(\uTj -\uTi)\zTi m_{ij}
-(\bb\cdot\nabla \uT,\zT)\rule{10ex}{0ex}\\
&&+\sum_{i\in\overline\Lambda} \int_{\Omega_i}[c_i\uTi \zTi -c\uT\zT]dx%\\
\end{eqnarray*}
\begin{eqnarray*}
&=&\sum_{i\in\overline\Lambda} \left\{\int_{\Omega_i}f(\zT-\zTi )dx+\int_{\Omega_i}(f-f_i)\zTi dx
\right\}\\
&&+\sumi \zTi \sumj\mu_{ij}(\uTi-\uTj )\dfrac{m_{ij}}{d_{ij}}-(A\nabla \uT, \nabla\zT)\\
&&+\sumi \Big\{\sumj (1-r_{ij})\gamma_{ij}(\uTj -\uTi)
\zTi m_{ij}\ -\ \int_{\Omega_i}(\bb\cdot\nabla \uT)\zTi dx\Big\}\\
&&-\sum_{i\in\overline\Lambda} \int_{\Omega_i}(\bb\cdot\nabla \uT)(\zT-\zTi )dx\\
&&+\sum_{i\in\overline\Lambda} \left\{\int_{\Omega_i}[c_i\uTi-c\uT]\zTi dx
\ -\ \int_{\Omega_i}c\uT(\zT-\zTi )dx\right\}\\
&=&\delta_0+\delta_1+\delta_2+\delta_3
\end{eqnarray*}
with
\begin{eqnarray*}
\delta_0&:=&\sumi \zTi \sumj\mu_{ij}(\uTi-\uTj )\dfrac{m_{ij}}{d_{ij}}
-(A\nabla \uT, \nabla\zT),\\
\delta_1&:=&\sum_{i\in\overline\Lambda} \int_{\Omega_i}
[f-\bb\cdot\nabla \uT-c\uT](\zT-\zTi )dx,\\
\delta_2&:=&\sumi \zTi \Big\{\int_{\Omega_i}[f-f_i+(\nabla\cdot\bb-c)\uT
+c_i\uTi]dx\ -\ \sumj \uTi\gamma_{ij}m_{ij}\Big\},\\
\delta_3&:=&\sumi \sumj \int_{\Gamma_{ij}}[(r_{ij}\uTi
+(1-r_{ij})\uTj )\gamma_{ij}-(\nu_{ij}\cdot\bb)\uT]\zTi ds.
\end{eqnarray*}
Here we have used that $\bb\cdot\nabla\uT=\nabla\cdot(\bb\uT)-(\nabla\cdot\bb)\uT.$

In the case of Donald diagrams, $\delta_0=0.$

In order to treat $\delta_0$ in the case of Voronoi diagrams,
we introduce a piecewise constant
(w.r.t.\ $\mathcal{T}$) approximation $A_\mathcal{T}$ to $A$ by
$\displaystyle A_\mathcal{T}|_T:=\frac{1}{\measd{T}}\int_T A\,dx,$
$T\in\mathcal{T}.$
Then we can write
\begin{eqnarray*}
\delta_0
&=&\sumi \zTi \sumj\left(\mu_{ij}-\frac{1}{m_{ij}}\int_{\Gamma_{ij}}A_\mathcal{T}ds\right)
(\uTi-\uTj )\dfrac{m_{ij}}{d_{ij}}\\
&&+\ \sumi \zTi \sumj\left(\int_{\Gamma_{ij}}A_\mathcal{T}ds\right)\dfrac{\uTi-\uTj}{d_{ij}}
-(A\nabla \uT, \nabla\zT).
\end{eqnarray*}
It is wellknown that, for arbitrary $\uT, \zT\in \VT,$
$$
\sumi \zTi \sumj\left(\int_{\Gamma_{ij}}A_\mathcal{T}ds\right)\dfrac{\uTi-\uTj}{d_{ij}}
=(A_\mathcal{T}\nabla \uT, \nabla\zT).
$$
Hence
$$
\delta_0
=\sumi \zTi \sumj\left(\mu_{ij}-\frac{1}{m_{ij}}\int_{\Gamma_{ij}}A_\mathcal{T}ds\right)
(\uTi-\uTj )\dfrac{m_{ij}}{d_{ij}}
+((A_\mathcal{T}-A)\nabla \uT, \nabla\zT).
$$
Since both $\nabla\uT, \nabla\zT$ are constant on every element $T\in\mathcal{T},$
the second term vanishes.
By a symmetry argument, we arrive at
$$
\delta_0
=\frac{1}{2}\sumi\sumj\left(\mu_{ij}-\frac{1}{m_{ij}}\int_{\Gamma_{ij}}A_\mathcal{T}ds\right)
(\uTi-\uTj)(\zTi-\zTj)\dfrac{m_{ij}}{d_{ij}}.
$$
Now the Cauchy-Schwarz inequality implies
\begin{eqnarray*}
|\delta_0|
&\le&\frac{1}{2}\left\{
\sumi\sumj\left(\mu_{ij}-\frac{1}{m_{ij}}\int_{\Gamma_{ij}}A_\mathcal{T}ds\right)^2
(\uTi-\uTj)^2\dfrac{m_{ij}}{d_{ij}}
\right\}^{1/2}\\
&&\times\left\{
\sumi\sumj(\zTi-\zTj)^2\dfrac{m_{ij}}{d_{ij}}
\right\}^{1/2}.
\end{eqnarray*}
The last factor can be bounded by $C_1|\zT|_{1,2,\Omega},$
therefore we get
\begin{equation}\label{eq:delta0est}
|\delta_0|\le C_1\eta_0|\zT|_{1,2,\Omega},
\end{equation}
where
$$
\eta_0^2:=\sumi\eta_{0i}^2
\quad\mbox{with}\quad
\eta_{0i}^2
:=\frac{1}{4}\sumj\left(\mu_{ij}-\frac{1}{m_{ij}}\int_{\Gamma_{ij}}A_\mathcal{T}ds\right)^2
(\uTi-\uTj)^2\dfrac{m_{ij}}{d_{ij}}\,.
$$
Setting $g:=f-\bb\cdot\nabla \uT-c\uT$ and
$\delta_{1i}:=\int_{\Omega_i}g(\zT-\zTi )dx,$
we can write (cf.\ Figure \ref{fig0} for the case $d=2$):
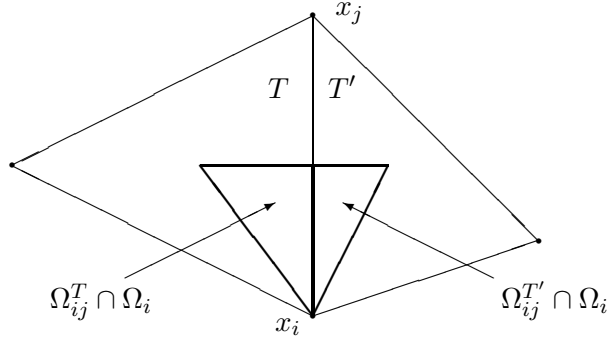
\begin{figure}
\begin{center}
\setlength{\unitlength}{1mm}
\begin{picture}(100,45)(-50,-2)
\thicklines
\put(0,0){\line(0,1){20}}
\put(0,0){\line(1,2){10}}
\put(0,0){\line(-3,4){15}}
\put(-15,20){\line(1,0){25}}
\thinlines
\put(0,20){\line(0,1){20}}
\put(0,0){\line(3,1){30}}
\put(0,0){\line(-2,1){40}}
\put(0,40){\line(-2,-1){40}}
\put(0,40){\line(1,-1){30}}
% nodes
\put(0,0){\circle*{0.7}}
\put(0,40){\circle*{0.7}}
\put(-40,20){\circle*{0.7}}
\put(30,10){\circle*{0.7}}
% notations
\put(-5,-2){$x_i$}
\put(3,40){$x_j$}
\put(-6,29){$T$}
\put(2,29){$T'$}
\put(-25,5){\vector(2,1){20}}
\put(-35,1){$\Omega_{ij}^T\cap\Omega_i$}
\put(24,5){\vector(-2,1){20}}
\put(25,1){$\Omega_{ij}^{T'}\cap\Omega_i$}
\end{picture}
\caption{The auxiliary simplices in the case $d=2$ for the Voronoi diagram}\label{fig0}
\end{center}
\end{figure}
$$
\delta_{1i}=\sumj \sum_{T\in\mathcal{T}\,:\,m_{ij}^T>0}
\int_{\Omega_{ij}^T\cap\Omega_i}g(\zT-\zTi )dx.
$$
On each simplex $T,$ it holds
$$
\zT=\zTi + \nabla\zT\cdot (x-x_i),
$$
where $ \nabla\zT$ is constant on $\Omega_{ij}^T.$

It follows that
\begin{eqnarray*}
\delta_{1i}&=&\sumj \sum_{T\in\mathcal{T}\,:\,m_{ij}^T>0}\int_{\Omega_{ij}^T\cap\Omega_i}g \nabla\zT\cdot
(x-x_i)dx\\
&\le&\sumj \sum_{T\in\mathcal{T}\,:\,m_{ij}^T>0}\int_{\Omega_{ij}^T\cap\Omega_i}| g|
\| \nabla\zT\|\|x-x_i\|dx\\
&\le&{\left\{\sumj \sum_{T\in\mathcal{T}\,:\,m_{ij}^T>0}\int_{\Omega_{ij}^T\cap\Omega_i}{| g|}^2
\|x-x_i\|^2dx\right\}}^{1/2}\\
&&\times
{\left\{\sumj \sum_{T\in\mathcal{T}\,:\,m_{ij}^T>0}\int_{\Omega_{ij}^T\cap\Omega_i}\| \nabla\zT\|^2dx
\right\}}^{1/2}\\
&\le&{\left\{\sumj \sum_{T\in\mathcal{T}\,:\,m_{ij}^T>0}h_T^2 %\|x_T^V-x_i\|^2
\int_{\Omega_{ij}^T\cap\Omega_i}{| g|}^2dx\right\}}
^{1/2}| \zT|_{1,2,\Omega_i}\\
&\le&\eta_{1i}| \zT|_{1,2,\Omega_i} ,
\end{eqnarray*}
where
$$
\eta_1^2:=\sum_{i\in\overline\Lambda}\eta_{1i}^2
\quad\mbox{with}\quad
\eta_{1i}^2:=\sumj \sum_{T\in\mathcal{T}\,:\,m_{ij}^T>0}h_T^2 %\|x_T^V-x_i\|^2
\int_{\Omega_{ij}^T\cap\Omega_i}g^2dx.
$$
Thus we arrive at
\begin{equation}\label{eq:8}
\delta_1\le\eta_1| \zT|_{1,2,\Omega}.
\end{equation}

For the third term $\delta_2,$ with
$$
\theta_i:=\int_{\Omega_i}[f-f_i+(\nabla\cdot\bb-c)\uT+c_i\uTi]dx\,-\,
\sumj \uTi\gamma_{ij}m_{ij},
$$
we have
$$
\delta_2=\sumi \zTi \theta_i.
$$
Because of
$$
\zTi \theta_i\le\eta_{2i}| \zTi |\sqrt{m_i},
$$
where $\eta_{2i}:=| \theta_i|/\sqrt{m_i},$ it follows
with
$\eta_2^2:=\sumi\eta_{2i}^2$
that
$$
\delta_2\le\eta_2\|\zT\|_\mathcal{T}.
$$
In view of the equivalence of the $L_2$-norm and the lumped
$L_2$-norm on $\VT,$ we obtain
\begin{equation}\label{eq:9}
\delta_2\le C_2\eta_2\|\zT\|_{0,2,\Omega}.
\end{equation}
For the remaining term $\delta_3$ we have, by a symmetry argument, that
$$
\delta_3=\sumi \delta_{3i},
$$
where
$$
\delta_{3i}:=\frac{1}{2}\sumj \int_{\Gamma_{ij}}\zetaTij (\zTi -\zTj)ds
$$
with
$$
\zetaTij :=[r_{ij}\uTi+(1-r_{ij})\uTj ]\gamma_{ij}\,-\,(\nu_{ij}\cdot\bb)\uT.
$$
In view of $\zTi -\zTj=d_{ij}(\nu_{ij}\cdot \nabla\zT)$ on $\Omega_{ij}^T$ we get
$$
\delta_{3i}=\frac{1}{2}\sumj d_{ij}\sum_{T\in\mathcal{T}\,:\,m_{ij}^T>0}
\int_{\Gamma_{ij}^T}\zetaTij (\nu_{ij}\cdot \nabla\zT)ds.
$$
It follows (remember that $\nu_{ij}\cdot \nabla\zT$ is constant on
$\Gamma_{ij}^T$ and
$ \nabla\zT$ is constant on $\Omega_{ij}^T\cap\Omega_i$)
\begin{eqnarray*}
\delta_{3i}&\le&\frac{1}{2}\sumj d_{ij}\sum_{T\in\mathcal{T}\,:\,m_{ij}^T>0}
\left|\int_{\Gamma_{ij}^T}\zetaTij ds\right|\| \nabla\zT\|\\
&=&\frac{1}{2}\sumj \sum_{T\in\mathcal{T}\,:\,m_{ij}^T>0}
\frac{d_{ij}}{\sqrt{\measd{\Omega_{ij}^T\cap\Omega_i}}}
\left|\int_{\Gamma_{ij}^T}\zetaTij ds\right|
\| \nabla\zT\|\sqrt{\measd{\Omega_{ij}^T\cap\Omega_i}}.
\end{eqnarray*}
By Cauchy's inequality, we have
\begin{eqnarray*}
\delta_{3i}&\le&\frac{1}{2}
{\left\{\sumj \sum_{T\in\mathcal{T}\,:\,m_{ij}^T>0}\frac{d_{ij}^2}{\measd{\Omega_{ij}^T\cap\Omega_i}}
{\left(\int_{\Gamma_{ij}^T}\zetaTij ds\right)}^2\right\}}^{1/2}| \zT|_{1,2,\Omega_i}\\
&\le&\eta_{3i}| \zT|_{1,2,\Omega_i},
\end{eqnarray*}
where
\begin{eqnarray*}
\eta_{3i}^2&:=&\frac{1}{4}\sumj
\sum_{T\in\mathcal{T}\,:\,m_{ij}^T>0}\frac{d_{ij}^2}{\measd{\Omega_{ij}^T\cap\Omega_i}}
{\left(\int_{\Gamma_{ij}^T}\zetaTij ds\right)}^2\\
&=&\frac{d}{4}\sumj \sum_{T\in\mathcal{T}\,:\,m_{ij}^T>0}\frac{d_{ij}}{m_{ij}^T}
{\left(\int_{\Gamma_{ij}^T}\zetaTij ds\right)}^2.
\end{eqnarray*}
Thus it holds that
\begin{equation}\label{eq:10}
\delta_3\le\eta_3| \zT|_{1,2,\Omega}.
\end{equation}
Summarizing the estimates (\ref{eq:delta0est}) -- (\ref{eq:10}), we obtain
$$
\eta_{nc}\le (C_1\eta_0+\eta_1+\eta_3)| \zT|_{1,2,\Omega}
+C_2\eta_2\| \zT\|_{0,2,\Omega},
$$
where the indicators have the following structure:
$$
\eta_l={\left\{\sumi \eta_{li}^2\right\}}^{1/2},
\qquad l\in\{0,1,2,3\},
$$
where
\begin{eqnarray*}
\eta_{0i}&=&\frac{1}{2}{\left\{\sumj\left(\mu_{ij}-\frac{1}{m_{ij}}
\int_{\Gamma_{ij}}A_\mathcal{T}ds\right)^2
(\uTi-\uTj)^2\dfrac{m_{ij}}{d_{ij}}\right\}}^{1/2}\\
&&\mbox{in case of Voronoi diagrams and $\eta_{0i}=0$ in case of Donald diagrams,}\\
\eta_{1i}&=&{\left\{\sumj \sum_{T\in\mathcal{T}\,:\,m_{ij}^T>0}h_T^2
\int_{\Omega_{ij}^T\cap\Omega_i}{[f-\bb\cdot\nabla \uT-c\uT]}^2dx
\right\}}^{1/2},\\
\eta_{2i}&=&\frac{1}{\sqrt{m_i}}\left|\int_{\Omega_i}[f-f_i+(\nabla\cdot\bb-c)
\uT+c_i\uTi]dx\,-\,\sumj \uTi\gamma_{ij}m_{ij}\right| ,\\
\eta_{3i}&=&{\left\{\frac{d}{2}\sumj \sum_{T\in\mathcal{T}\,:\,m_{ij}^T>0}
\frac{d_{ij}}{m_{ij}^T}
{\left(\int_{\Gamma_{ij}^T}[(r_{ij}\uTi+(1-r_{ij})\uTj )\gamma_{ij}-
(\nu_{ij}\cdot\bb)\uT]ds\right)}^2\right\}}^{1/2}.
\end{eqnarray*}
\begin{remark}
(i) We mention that all the indicators $\eta_l$ can be rewritten in such a way
that the resulting local indicators are related to the elements $T\in\mathcal{T}.$

(ii) It can be shown that the indicators $\eta_l$ are \emph{order-consistent}
with the \emph{a priori} error estimate (Theorem \ref{errortheorem})
in the following sense :

If $f\in W^1_q(\Omega)$ with some $q>d$ and $u\in W^2_2(\Omega),$
then there is a constant $C_c>0$ such that
$$
\sum_{l=0}^3 \eta_l\le C_ch\left[\|u\|_{2,2}+\|f\|_{1,r}\right],
$$
see \cite[Thm.\ 4]{Angermann:92a} for a special case.
\end{remark}

\section{Conclusions and perspectives}

We derived an estimator for measuring simultaneously two types of errors,
modeling and discretization errors, with respect
to user-defined output functionals.
The approach is formulated for stationary nonlinear partial differential equations
involving complex models.
The main focus was on the consideration of discretization methods
which do not possess the property of Galerkin orthogonality.
For the example of node-centered finite volume methods,
by localization of the estimators we presented local error indicators
which allow for local mesh refinement and local model modification.

In future work, the described framework will be the basis
for a more detailed analysis of problems with tensor-valued diffusion coefficients
and dominating convection and for algorithms which balance the
indicators corresponding to the different sources of error.

\newcommand{\etalchar}[1]{$^{#1}$}

\end{document}